\documentclass[11pt,notitlepage,twoside]{article}
\usepackage{latexsym}
\usepackage{amsmath}
\usepackage{amsfonts}
\usepackage{amssymb}

\newcommand{\be}{\begin{equation}}
\newcommand{\ee}{\end{equation}}
\newcommand{\ba}{\begin{align*}}
\newcommand{\ea}{\end{align*}}

\newcommand{\R}{\mathbb{R}}
 
\newtheorem{theorem}{Theorem}[section]

\newtheorem{proposition}{Proposition}[section]

\author{Zakaria Boucheche\\
{\footnotesize  Laboratory of Applied Mathematics and Harmonic
Analysis (LR11-ES52)}\\{\footnotesize  (Faculty of Sciences of Gabes.
Gabes University)}}

\title { \Large \textbf{An improved existence criterion and an
optimal result}}

\begin{document}

\date{ }

\maketitle

\begin{abstract}We are concerned with a semi-linear elliptic equation on a smooth bounded domain $\Omega$ of $\mathbb{R}^n,\,n\geq
5,$ which involves a critical nonlinearity and a linear term of the
form $K(x)u^{(n+2)/(n-2)}$ and $\mu u,$ respectively. By using a
test function procedure, we give an existence criterion involving
the parameter $\mu$ and the function $K(x).$ For a particular case
of $\Omega,\,K(x)$ and $n,$ we prove its optimality through a
Pohozaev type identity.
\end{abstract}

\noindent{ \footnotesize {{\it \\\\2010 Mathematics Subject
Classification.} $\quad$ 35J60, 58K05, 35J08, 41A58, 49K15.}
}
\section{Introduction and main results}
\def\theequation{1.\arabic{equation}}\makeatother
\setcounter{equation}{0}
 This paper is a second part devoted to the study of the following nonlinear elliptic partial differential
equation with zero Dirichlet boundary condition
\begin{equation}\label{problem}
\begin{aligned}
-\Delta& u=K(x)u^{q}+\mu u \quad \mbox{in}\,\,\Omega,\\{}& u>0 \quad
\mbox{in} \,\,\Omega,\quad u=0 \quad \mbox{on} \,\,\partial\Omega,
\end{aligned}
\end{equation}where $\Omega\subset \mathbb{R}^n,\,n\geq 5,$ is a bounded domain with a smooth boundary
$\partial\Omega,$ $K(x)$ is a $C^2$-function in $\bar{\Omega},$
$q+1=\frac{2n}{n-2}$ is the critical exponent for the embedding
$H_0^1\bigl(\Omega\bigr)$ into $L^{q+1}\bigl(\Omega\bigr)$ and $0<
\mu< \mu_1(\Omega),$ where $\mu_1(\Omega)$ denotes the first
eigenvalue of $(-\Delta)$ in
$H_0^1\bigl(\Omega\bigr),$ \\

In \cite[Theorem 1.1]{Bo1}, we were interested on the existence of
at least one solution to \eqref{problem}. This result was centered
on a Lions's condition. Namely, by using, we have proved the
following theorem. Denote $\sup_{\bar{\Omega}}(K(x)):=K_\infty$ and
$S:=\mathrm{inf}\{\|u\|^2,\,\,\, u\in
H_0^1\bigl(\Omega\bigr)\,\,\mathrm{and\,\,}\|u\|_{q+1}=1\}$ is the
best Sobolev constant, where $J(u):=\int_\Omega
K(x)|u(x)|^{q+1}\,\mathrm{d}x$ and $\|u\|_p^p:=\int_\Omega
|u(x)|^p\,\mathrm{d}x$ for any $p> 1.$
\begin{theorem}\label{Lions}$\mathrm{\bigl([4]\bigr)}$
\begin{equation}\label{zakaria}
\quad
\end{equation}\end{theorem}

When $K(x)\equiv 1,$ we recognize the Brezis--Nirenberg existence result
\cite[Lemma 1.2]{BrNir}. In order to establish the condition
\eqref{zakaria}, Brezis and Nirenberg \cite[Lemma 1.1]{BrNir} follow
an original idea due to Aubin \cite{A}: By considering the following
test function
\begin{equation}\label{aubin}u_{\lambda,y_0}(x)=\varphi(x)\cdot c_n^{\frac{n-2}{4}}
\bigl(\frac{\lambda}{1+\lambda^2|x-y_0|^2}\bigr)^\frac{n-2}{2}=:\varphi(x)\cdot
\delta_{y_0,\lambda}(x),\end{equation}where $c_n:=n^2-2n,\,y_0\in
\Omega,\,\lambda> 0,\,\delta_{y_0,\lambda}$ are the positive
solutions in $ \mathbb{R}^n,$ concentrated at $y_0,$ of $-\Delta u =
u^{\frac{n+2}{n-2}}$ and $\varphi$ is a cut-off function, they proved that the condition \eqref{zakaria} is satisfied for any $\mu> 0.$\\

When $K(x)\not\equiv 1,$ the situation becomes extremely different:
Indeed the behavior of $K(x)$ plays a crucial role in establishing
existence results; see, e.g., \cite{Bo} for the case $\mu=0$ and
$K(x)$ is positive everywhere. But for $\mu=0$ and, of course,
$K_\infty>0,$ the following Pohozaev identity \cite{P}
\begin{equation}\label{Pohozaev}\frac{1}{2}\int_{\partial
\Omega}|\frac{\partial u}{\partial\nu}(x)|^2\langle x,\,
\nu(x)\rangle\,\mathrm{d}x =\frac{n-2}{2n}\int_{\Omega}\bigl\langle
x,\, \nabla K(x)\bigr\rangle
u^{\frac{2n}{n-2}}(x)\,\mathrm{d}x+\mu\int_{\Omega}u^2(x)\,\mathrm{d}x\end{equation}
asserts that the problem \eqref{problem} has no solution provided
that $\Omega$ is star-shaped with respect to the origin $o$ of
$\mathbb{R}^n$ and $\langle x,\, \nabla K(x)\rangle\leq 0$ in
$\Omega.$ Here $\nu(x)$ denotes the outward normal vector at $x$ to
$\partial \Omega$ and $u$ is supposed to be a solution of
\eqref{problem}. (This identity \eqref{Pohozaev} is obtained by
multiplying the equation given in \eqref{problem} on the one hand by
$u$ and on the other hand by $\sum_{i=1}^nx_i(\partial u/\partial
x_i$), and using an integration by parts and the fact that on
$\partial \Omega$ we have $\nabla u=(\partial u/\partial \nu)\nu$).
In fact, in this case the condition \eqref{zakaria} is not
satisfied. In view of this nonexistence result, naturally one can
ask: What is the concrete condition that can we impose on $\mu$ and
an absolute maximum $y_0$ of $K(x)$ in $\bar{\Omega}$ so that
\eqref{zakaria} becomes satisfied? In the case $\mu> 0,$ Lions
\cite[Remark 4.7]{Li} considered the test function \eqref{aubin} and
he showed that the condition \eqref{zakaria} is satisfied provided
that
\begin{eqnarray*}
K_\infty=K(y_0)>
0\quad \mathrm{with}\quad y_0 \in \Omega,\,\,\,\quad\qquad\label{lionss}\\
-\frac{(n-2)^2\bar{c}_2\Delta K(y_0)}{2nK(y_0)}< \mu
\bar{c}_3,\quad\qquad\qquad\label{lions}
\end{eqnarray*}
where $\bar{c}_2$ and $\bar{c}_3$ are two positive constants
depending only on $n;$ see Proposition \ref{proposition} below.\\

Convinced to expand the validity of the condition \eqref{zakaria} to
more large class of functions $K(x)$ when $\mu$ is fixed, a choice
of a test function taking care of the geometry of $\Omega$ becomes
useful.  To this end, let $P$ be the projection from $H^1(\Omega)$
onto $H_0^1(\Omega)$; that is, $u= Pf$ is the unique solution of
$\Delta u = \Delta f \,\,\mbox{in}\,\,\Omega,\,\,\, u=0\,\,
\mbox{on}\,\,\partial\Omega.$ Denote by $H$ the regular part of the
Green's function of $\bigl(-\Delta\bigr)$ on $\Omega.$ By using the
test function $P\delta_{y_0,\lambda},$ we are able to prove the
following proposition:
\begin{proposition}\label{proposition}{\it\,\,Let $n\geq 5.$ Let $K(x)\in C^2(\bar{\Omega})$ satisfying $K_\infty=K(y_0)>
0$ with $y_0\in \Omega$ and let $\mu> 0.$ Then the condition
\eqref{zakaria} holds true provided that one of
the following two conditions is satisfied:\\
${\bf i)}$  $\quad-\frac{(n-2)^2\bar{c}_2\Delta K(y_0)}{2nK(y_0)}<
\mu \bar{c}_3,$\\\\ ${\bf ii)}$ $\quad-\frac{(n-2)^2\bar{c}_2\Delta
K(y_0)}{2nK(y_0)}= \mu \bar{c}_3$ and
\begin{equation*}\label{not satisfied}
\begin{aligned}&\liminf_{\lambda \rightarrow +\infty}\,\lambda^{n-2} \biggl[-\int_{B_0}\bigl(\frac{K(x)}{K(y_0)}-1-\frac{\Delta
K(y_0)}{2nK(y_0)}|x|^2\bigr)\delta_{y_0,
\lambda}^{\frac{2n}{n-2}}\mathrm{d}x\Bigr]+S_n\sum_{k=2}^{[\frac{n-2}{2}]}a_{n,k}\frac{\mu^k}{\lambda^{2k}}\biggr]\\&<
\frac{n\bar{c}_4}{(n-2)},\end{aligned}\end{equation*}where
$d_0:=\mathrm{dist}(y_0,\partial \Omega),$ $B_0$ is the ball of
center $y_0$ and radius $d_0,$
$S_n:=\int_{\mathbb{R}^n}\bigl(1+|x|^2\bigr)^{-n}\mathrm{d}x
,\,\,\bar{c}_2=\int_{\mathbb{R}^n} |x|^2/(1+|x|^2)^n
\mathrm{d}x,\,\,\bar{c}_3=\int_{\R^n}1/(1+|x|^2)^{n-2}\,\mathrm{d}x,$
$a_{n,k}$'s are the constants defined by the following Taylor
expansion
$$\bigl(1-\frac{\bar{c}_3}{c_nS_n}t\bigr)^{\frac{n}{n-2}}
=1-\frac{n\bar{c}_3}{(n-2)c_nS_n}t+\sum_{k=2}^{[\frac{n-2}{2}]}a_{n,k}t^k+o\bigl(t^{\frac{n-2}{2}}\bigr)\quad
as\quad t\rightarrow 0,$$\begin{equation*}\label{specified}
\begin{aligned}\bar{c}_4:=&-H(y_0,
y_0)\int_{\mathbb{R}^n}
\frac{\mathrm{d}x}{(1+|x|^2)^{\frac{n+2}{2}}}+\mu
c_n^{-1}\biggl[2\int_{\Omega}H(y_0,\,x)\frac{\mathrm{d}x}{|x-y_0|^{n-2}}\\&-
\int_{\Omega}H^2(y_0,\,x)\mathrm{d}x+\int_{\R^n\setminus
\Omega}\frac{\mathrm{d}x}{|x-y_0|^{2n-4}}\biggr]\end{aligned}\end{equation*}
}
\end{proposition}
{\bf Remark 1.1}\,\,\,\,If we use the test function
$u_{\lambda,y_0}$ instead of $P\delta_{y_0,\lambda},$ then the
corresponding constant $\bar{c}_4$ can not be specified. This is due
to the fact that
$u_{\lambda,y_0}$ does not deal with the boundary of $\Omega.$\\\\
{\bf Example 1.1}\,\,\,\,Let $\mu> 0$ be fixed. To verify if a
function $K(x)$ satisfies the condition ${\bf (ii)},$ we need to
know its Taylor expansion, near $y_0,$ of order greater than $2.$
For example, let us take the case $\Omega=B$ is the unit ball of
$\mathbb{R}^5$ and $y_0$ is the origin of $\mathbb{R}^5.$ Assume
that $K(x)=f(|x|)$ is radial and radially non-increasing function
with$$f(t)=f(0)+at^2+bt^{3}+o(t^3)\quad \mathrm{as}\,\,t\rightarrow
0.$$Then the condition ${\bf (ii)}$ is satisfied provided
that$$-9\bar{c}_2a= \mu \bar{c}_3f(0)\quad \mathrm{and}\quad
-3b\int_{\mathbb{R}^5} |x|^3/(1+|x|^2)^5 \mathrm{d}x<
5\bar{c}_4f(0).$$\\

In a second part of this work, we will try to analyze the optimality
of the condition ${\bf (i)}$ for some class of functions when
$\Omega$ is a ball, $n$ is odd with $n=5$ or $n> 19$ and $K(x)$ is
close to a constant,
radial and radially non-increasing. To this end, let us state the following assumptions: Assume that\\

$\mathbf{(K_1)}$  $\Omega=B(y_0,\,\gamma)$ is the ball of center $y_0$ and radius $\gamma$ in $\mathbb{R}^n.$\\

$\mathbf{(K_\eta)}$  $K(x)=K(y_0)+\eta f_1(|x-y_0|)$ is a
non-negative $C^2$-function in $\bar{\Omega},$ where $\eta,\,K(y_0)>
0$ are fixed
constants and $f_1$ is a non-increasing function on $[0,\,\gamma]$ independent of
$\eta.$\\\\
In this case, we will refer to the problem \eqref{problem}
as $\mathbf{(BN)}_\eta.$\\

$\mathbf{(K_3)}$  $\limsup\limits_{t\rightarrow
0}\frac{f_1'(t)-f_1''(0)t}{t^{n-3}}< +\infty.$\\\\
Our optimal result is the following:
\begin{theorem}\label{nonexistence}{\it\,\,Let $n$ be an odd integer with $n=5$ or $n> 19$ and let $0<
\mu< \mu_1(\Omega).$ Assume that $\Omega$ and $K(x)$ satisfy the
assumptions $\mathbf{(K_1)},\,\mathbf{(K_\eta)}$ and
$\mathbf{(K_3)}.$ Then there exists a constant $\bar{\eta}$
depending on $n,\,f_1(t)$ and $K(y_0)$ such that if $0< \eta\leq
\bar{\eta},$ then the problem $\mathbf{(BN)}_\eta$ admits a solution
if and only
if\begin{equation}\label{lions}-\frac{(n-2)^2\bar{c}_2\Delta
K(y_0)}{2nK(y_0)}< \mu \bar{c}_3.\end{equation}}\end{theorem} The
proof of the sufficiency is obtained by a combination of the results
Theorem \ref{Lions} and Proposition \ref{proposition}. For the
necessity of the condition \eqref{lions}, we argue by contradiction:
The key point is to establish an adequate Pohozaev type identity for
the desired solution of \eqref{problem}; this identity is a natural
extension to that given in the proof of \cite[Lemma 1.4]{BrNir}. To
conclude, we need to investigate a constant $\bar{\eta}$ depending
only on $n,\,f_1(t)$ and $K(y_0)$ such that if $\eta\leq \bar{\eta}$
and $\bigl[-(n-2)^2\bar{c}_2\Delta K(y_0)\bigr]/2nK(y_0)\geq \mu
\bar{c}_3,$ then this identity becomes impossible.
\section{Proof of the results}\label{sec2}
\def\theequation{2.\arabic{equation}}\makeatother
\setcounter{equation}{0} {\it Proof of Proposition
\ref{proposition}.} Let $y_0\in \Omega$ be such that
$K_\infty=K(y_0)> 0.$ Denoting, for $\lambda> 0$ a fixed constant
large
enough,\begin{equation}\label{Bahri}A_{y_0,\mu}(\lambda):=\frac{\int_{\Omega}|\nabla
P\delta_{y_0, \lambda}|^2- \mu \int_{\Omega} (P\delta_{y_0,
\lambda})^2} {\Bigl( \int_{\Omega}K(P\delta_{y_0,
\lambda})^{\frac{2n}{n-2}} \Bigr)^{\frac{n-2}{n}}}.
\end{equation}In order to get the claim of Proposition
\ref{proposition}, it is sufficient to prove that, for $\lambda$
large enough,
\begin{equation}\label{Bahri1}\bigl[A_{y_0,\mu}(\lambda)\bigr]^{\frac{n}{n-2}}<
\frac{1}{K(y_0)}S^{\frac{n}{n-2}}
\end{equation}provided that one of the conditions ${\bf (i)}$ and ${\bf (ii)}$ is satisfied. To this end, we need an estimation of the following
three quantities:$$\int_{\Omega} (P\delta_{y_0, \lambda})^2,\quad
\int_{\Omega}K(P\delta_{y_0, \lambda})^{\frac{2n}{n-2}} \quad
\mathrm{and}\quad \int_{\Omega}|\nabla P\delta_{y_0, \lambda}|^2.$$
The last two quantities were estimated in \cite[(2.67), (5.31) and
Estimate F8]{B}, and we have\begin{eqnarray}
\int_{\Omega}K(P\delta_{y_0,
\lambda})^{\frac{2n}{n-2}}=K(y_0)c_n^{\frac{n}{2}}\Bigl[S_n
+\frac{\bar{c}_2}{2n}\frac{\Delta K(a)}{K(y_0)\lambda^2}-
\frac{2n\bar{c}_1}{n-2}\frac{H(y_0, y_0)}{\lambda^{n-2}}\nonumber\\
\qquad\quad\qquad\qquad+\int_{B_0}\bigl(\frac{K(x)}{K(y_0)}-1-\frac{\Delta
K(y_0)}{2nK(y_0)}|x|^2\bigr)\delta_{y_0,
\lambda}^{\frac{2n}{n-2}}\mathrm{d}x\Bigr]\label{Bahri2}\\+o(\frac{1}{\lambda^{n-2}})+
O\Bigl(\frac{log(\lambda d_0)}{(\lambda d_0)^n}\Bigr),\qquad\quad\quad\qquad\qquad\nonumber\\
\int_{\Omega}|\nabla P\delta_{y_0,
\lambda}|^2=c_n^{\frac{n}{2}}\bigl[S_n - \bar{c}_1\frac{H(y_0,
y_0)}{\lambda^{n-2}}\bigr] + O\Bigl(\frac{log(\lambda d_0)}{(\lambda
d_0)^n}\Bigr),\quad\qquad\,\,\,\label{Bahri2222}
\end{eqnarray}
where $d_0:=\mathrm{dist}(y_0,\,\partial \Omega),\,B_0$ is the ball
of center $y_0$ and radius $d_0,$ $\bar{c}_1=\int_{\mathbb{R}^n}
\frac{dx}{(1+|x|^2)^{\frac{n+2}{2}}}$ and
$\bar{c}_2=\int_{\mathbb{R}^n} \frac{|x|^2}{(1+|x|^2)^n} dx.$ Then
we are left with the first
quantity:\begin{equation}\label{Bahri3}\int_\Omega(\mathrm{P}\delta_{y_0,\lambda})^2(x)\,\mathrm{d}x=
\int_\Omega\delta^2_{y_0,\lambda}\,\mathrm{d}x+\int_\Omega\theta^2_{y_0,\lambda}(x)\,\mathrm{d}x-2\int_\Omega\delta_{y_0,\lambda}\theta_{y_0,\lambda}(x)\,\mathrm{d}x,\end{equation}
where $\,\theta_{y_0,\lambda}:= \delta_{y_0, \lambda} -
P\delta_{y_0, \lambda}.$ First we recall that from \cite[(5.25)]{B}
we have the following estimate
\begin{equation*}\label{Bahri4}\theta_{y_0,\lambda} (x)=
\frac{c_n^{\frac{n-2}{4}}}{\lambda^{\frac{n-2}{2}}}H(y_0, x) +
\frac{1}{\lambda^{\frac{n+2}{2}}d_0^n}\cdot
O(1),\quad\forall\,\,x\in \Omega,\end{equation*}where $|O(1)|$ is a
quantity upper-bounded by a positive constant $M$ independent of
$x\in \Omega.$ This, together with Lebesgue's dominated convergence
theorem, implies that
\begin{eqnarray}
\int_\Omega\delta_{y_0,\lambda}\theta_{y_0,\lambda}(x)\,\mathrm{d}x=
\frac{c_n^{\frac{n-2}{2}}}{\lambda^{n-2}}\int_{\Omega}H(y_0,\,x)\frac{\lambda^{n-2}}{(1+\lambda^2|x-y_0|^2)^{\frac{{n-2}}{2}}}\,\mathrm{d}x+o(
\frac{1}{\lambda^{n-2}})\nonumber\\=\frac{c_n^{\frac{n-2}{2}}}{\lambda^{n-2}}\int_{\Omega}H(y_0,\,x)\frac{1}{|x-y_0|^{n-2}}\,\mathrm{d}x+o(
\frac{1}{\lambda^{n-2}}),\,\,\,\quad\qquad\label{Bahri5}\\
\int_\Omega\theta^2_{y_0,\lambda}\mathrm{d}x=\frac{c_n^{\frac{n-2}{2}}}{\lambda^{n-2}}\int_{\Omega}H^2(y_0,\,x)\mathrm{d}x+o(
\frac{1}{\lambda^{n-2}}),\quad\qquad\qquad\label{Bahri6}
\end{eqnarray}On the other hand, by using, again, Lebesgue's dominated convergence theorem
we
obtain\begin{equation}\label{regu}\begin{aligned}\int_\Omega\delta^2_{y_0,\lambda}\,\mathrm{d}x
&=c_n^{\frac{n-2}{2}}\biggl(\frac{\bar{c}_3}{\lambda^2}-
\int_{\R^n\setminus
\Omega}\frac{\lambda^{n-2}}{(1+\lambda^2|x-y_0|^2)^{n-2}}\,\mathrm{d}x\biggr)\\&
=c_n^{\frac{n-2}{2}}\biggl(\frac{\bar{c}_3}{\lambda^2}-\frac{\bar{\bar{c}}_5}{\lambda^{n-2}}+o(
\frac{1}{\lambda^{n-2}})\biggr),\end{aligned}\end{equation} where
$\bar{c}_3=\int_{\R^n}\frac{1}{(1+|x|^2)^{n-2}}\,\mathrm{d}x$ and
$\bar{\bar{c}}_5=\int_{\R^n\setminus
\Omega}\frac{1}{|x-y_0|^{2n-4}}\,\mathrm{d}x.$ Combining
\eqref{Bahri3}--\eqref{regu} we obtain
\begin{eqnarray}
\int_\Omega(\mathrm{P}\delta_{y_0,\lambda})^2(x)\,\mathrm{d}x=
\frac{c_n^{\frac{n-2}{2}}}{\lambda^{n-2}}\biggl[-2\int_{\Omega}H(y_0,\,x)\frac{1}{|x-y_0|^{n-2}}\,\mathrm{d}x+
\int_{\Omega}H^2(y_0,\,x)\mathrm{d}x\biggr]\nonumber\\+
c_n^{\frac{n-2}{2}}\biggl(\frac{\bar{c}_3}{\lambda^2}-\frac{\bar{\bar{c}}_5}{\lambda^{n-2}}\biggr)+o(
\frac{1}{\lambda^{n-2}})\quad
\qquad\qquad\qquad\quad\quad\nonumber\\
=:c_n^{\frac{n-2}{2}}\biggl(\frac{\bar{c}_3}{\lambda^2}-\frac{\bar{c}_6}{\lambda^{n-2}}\biggr)+o(
\frac{1}{\lambda^{n-2}}),\quad
\qquad\qquad\qquad\quad\qquad\quad\label{Bahri7}
\end{eqnarray}where $\bar{c}_6:=\int_{\Omega}\bigl[2H(y_0,\,x)/|x-y_0|^{n-2}-
H^2(y_0,\,x)\bigr]\mathrm{d}x+\bar{\bar{c}}_5.$ Combining
\eqref{Bahri7} and \eqref{Bahri2222} we get, for $\lambda$ large
enough,
\begin{equation}\label{Bahri8}\begin{aligned}
&\biggl[\int_{\Omega}|\nabla P\delta_{y_0, \lambda}|^2-\mu
\int_\Omega(\mathrm{P}\delta_{y_0,\lambda})^2(x)\,\mathrm{d}x\biggr]^{\frac{n}{n-2}}\\&=(S_nc_n^{\frac{n}{2}})^{\frac{n}{n-2}}
\biggl[1-\frac{\mu
\bar{c}_3}{c_nS_n\lambda^2}-\frac{\bar{c}_0}{\lambda^{n-2}}+o(
\frac{1}{\lambda^{n-2}})\biggr]^{\frac{n}{n-2}}\\&=(S_nc_n^{\frac{n}{2}})^{\frac{n}{n-2}}\biggl[1-\frac{n\mu
\bar{c}_3}{(n-2)c_nS_n\lambda^2}-\frac{n\bar{c}_0}{(n-2)\lambda^{n-2}}+\sum_{k=2}^{[\frac{n-2}{2}]}a_{n,k}\frac{\mu^k}{\lambda^{2k}}\biggr]+o(
\frac{1}{\lambda^{n-2}}),
\end{aligned}\end{equation}
where $a_{n,k}$'s are fixed constants defined by the following
Taylor
expansion$$\bigl(1-\frac{\bar{c}_3}{c_nS_n}t\bigr)^{\frac{n}{n-2}}
=1-\frac{n\bar{c}_3}{(n-2)c_nS_n}t+\sum_{k=2}^{[\frac{n-2}{2}]}a_{n,k}t^k+o\bigl(t^{\frac{n-2}{2}}\bigr)\quad
as\quad t\rightarrow 0$$($[(n-2)/2]$ denotes the integer part of
$(n-2)/2$ and the sum $\sum_{k=2}^{[(n-2)/2]}$ is omitted when
$n=5$) and $\bar{c}_0:=\bigl(\bar{c}_1H(y_0, y_0)+\mu
c_n^{-1}\bar{c}_6\bigr)/S_n.$ \eqref{Bahri}, \eqref{Bahri2} and
\eqref{Bahri8} imply that, for $\lambda$ large enough,
\begin{equation}\label{Bahri9}\begin{aligned}
&\bigl[A_{y_0,\mu}(\lambda)\bigr]^{\frac{n}{n-2}}\\&=\frac{S_n^{\frac{2}{n-2}}c_n^{\frac{n}{n-2}}}{K(y_0)}\biggl[1-\bigl(\frac{\bar{c}_2}{2n}\frac{\Delta
K(a)}{K(y_0)}+\frac{\mu
\bar{c}_3}{(n-2)^2}\bigr)\frac{1}{S_n\lambda^2}-\frac{n\bar{c}_4}{(n-2)S_n\lambda^{n-2}}+o(
\frac{1}{\lambda^{n-2}})\\&\quad\quad\qquad\quad\quad-\frac{1}{S_n}\int_{B_0}\bigl(\frac{K(x)}{K(y_0)}-1-\frac{\Delta
K(y_0)}{2nK(y_0)}|x|^2\bigr)\delta_{y_0,
\lambda}^{\frac{2n}{n-2}}\mathrm{d}x\Bigr]+\sum_{k=2}^{[\frac{n-2}{2}]}a_{n,k}\frac{\mu^k}{\lambda^{2k}}\biggr],
\end{aligned}\end{equation}
where $\bar{c}_4:=-\bar{c}_1H(y_0, y_0)+\mu c_n^{-1}\bar{c}_6.$ On
the other hand, observe that since $K(x)\in C^2(\bar{\Omega}),$ then
\begin{equation}\label{Bahri999}\int_{B_0}\bigl(\frac{K(x)}{K(y_0)}-1-\frac{\Delta
K(y_0)}{2nK(y_0)}|x|^2\bigr)\delta_{y_0,
\lambda}^{\frac{2n}{n-2}}\mathrm{d}x=o(\frac{1}{\lambda^2}).\end{equation}Observe
also that
\begin{equation}\label{Bahri000}S=c_nS_n^{\frac{2}{n}}.\end{equation}Thus
under the condition ${\bf (i)}$, the claim \eqref{Bahri1} follows by
combining \eqref{Bahri9}--\eqref{Bahri000} and taking $\lambda$
large enough. If the condition ${\bf (ii)}$ is satisfied instead of
${\bf (i)},$ then \eqref{Bahri1} follows by taking $\lambda$ large
enough in the right hand side of \eqref{Bahri9} and using
\eqref{Bahri000}. This finishes the proof of Proposition
\ref{proposition}.\\
{\it Proof of Theorem \ref{nonexistence}.} {\it\, Sufficiency of the
condition \eqref{lions}:} From $\mathbf{(K_1)}$ and
$\mathbf{(K_\eta)}$ we get $K_\infty=K(y_0)> 0$ with $y_0\in
\Omega.$ This, together with the condition \eqref{lions} and the
result of Proposition \ref{proposition}, implies that
\eqref{zakaria} is
satisfied. Thus a solution to problem \eqref{problem} is obtained by applying Theorem \ref{Lions}.\\
{\it Necessity of the condition \eqref{lions}:}\,\,Arguing by
contradiction, assuming that the problem \eqref{problem} has a
solution $u$ under the
condition\begin{equation}\label{condition}-\frac{(n-2)^2\bar{c}_2\Delta
K(y_0)}{2nK(y_0)}\geq \mu \bar{c}_3.\end{equation}In particular, we
have $\Delta K(y_0)\neq 0.$ Up to a translation and a dilatation in
the space, we can suppose that
$$\Omega=B\quad \mathrm{is\, the\, unit\, ball\, of}\,\,
\mathbb{R}^n.$$ Now, by a result of Gidas--Ni--Nirenberg
\cite[Theorem 1$'$]{GNNir}, $\mathbf{(K_1)}$ and $\mathbf{(K_\eta)}$
imply that $u$ is necessarily spherically symmetric. We write
$u(x)=:u(t)$ and $K(x)=:f(t),$ where $t=|x|\in [0,\,1].$ Thus $u$
satisfies the following ordinary differential equation
\begin{eqnarray}
&\qquad\qquad-u{''}-\frac{n-1}{t}u'=f(t)u^{\frac{n+2}{n-2}}+\mu
u\qquad \mathrm{on\,\,}(0,\,1), \label{zak1}\\& u'(0)=u(1)=0.
\nonumber
\end{eqnarray}(Note that $u\in C^2(\bar{\Omega})$).  Let $\psi$ be a smooth function on $[0,\,1]$ such
that $\psi(0)=0.$ Multiplying the equation \eqref{zak1} by
$t^{n-1}\psi u'$ and
$\Bigl(t^{n-1}\psi'(t)-(n-1)t^{n-2}\psi(t)\Bigr)u$ and integrating
by parts several times in order to obtain \begin{eqnarray}
-\frac{1}{2}|u'(1)|^2\psi(1)+\frac{1}{2}\int_0^1|u'(t)|^2\Bigl(t^{n-1}\psi'(t)-(n-1)t^{n-2}\psi(t)\Bigr)\mathrm{d}t\nonumber\\
=-\bar{c}_n\int_0^1u^{\frac{2n}{n-2}}\biggl[f(t)\Bigl(t\psi'(t)+(n-1)\psi(t)\Bigr)+f'(t)t\psi(t)\biggr]t^{n-2}\mathrm{d}t
\label{zak3}\\-\frac{\mu}{2}
\int_0^1u^{2}\Bigl(t^{n-1}\psi'(t)+(n-1)t^{n-2}\psi(t)\Bigr)\mathrm{d}t
,\qquad\qquad\qquad\,\,\,\,\nonumber\\
\int_0^1\biggl[f(t)u^{\frac{2n}{n-2}}\Bigl(t\psi'(t)-(n-1)\psi(t)\Bigr)+\mu
u^{2}\Bigl(t\psi'(t)-(n-1)\psi(t)\Bigr)\biggr]t^{n-2}\mathrm{d}t\nonumber\\
=-\frac{1}{2}\int_0^1u^{2}\biggl[t^3\psi^{(3)}(t)+(n-1)(n-3)\Bigl(\psi(t)-t\psi'(t)\Bigr)\biggr]t^{n-4}\mathrm{d}t\label{zak4}\\
+\int_0^1|u'(t)|^2\Bigl(t\psi'(t)-(n-1)\psi(t)\Bigr)t^{n-2}\mathrm{d}t,\qquad\qquad\qquad\quad\,\,\,\,\nonumber
\end{eqnarray}respectively, where
$\bar{c}_n:=\frac{n-2}{2n}.$ Combining \eqref{zak3} and \eqref{zak4}
we get
\begin{align}\label{zak5}\begin{split}
&-\frac{1}{2}|u'(1)|^2\psi(1)+\int_0^1u^{2}\biggl[\mu\psi'(t)+\frac{1}{4}\psi^{(3)}(t)
+\frac{1}{4}\frac{(n-1)(n-3)}{t^3}\bigl(\psi(t)-t\psi'(t)\bigr)\biggr]t^{n-1}\mathrm{d}t\\&=
\int_0^1u^{\frac{2n}{n-2}}\biggl[-\bar{c}_ntf'(t)\psi(t)+\frac{(n-1)}{n}f(t)\Bigl(\psi(t)-r\psi'(t)\Bigr)\biggr]t^{n-2}\mathrm{d}t.
\end{split}\end{align}(Note that the Pohozaev identity \eqref{Pohozaev}
corresponds to the case where $\psi(t)=t$). In order to get the
desired contradiction, we need to choose a suitable function $\psi$
as a solution of the following ordinary differential
equation\begin{equation}\label{zak6}
\mu\psi'+\frac{1}{4}\psi^{(3)}+\frac{1}{4}\frac{(n-1)(n-3)}{t^3}\bigl(\psi-t\psi'\bigr)=0,\quad\forall\,\,t\in
(0,\,1].
\end{equation}
A straightforward computation shows that the equation \eqref{zak6}
has two solutions defined on $[0,\,1]$ by a series
$\psi_1(t)=\sum_{p=0}^{+\infty}a_{2p+1}t^{2p+1}$ and
$\psi_2(t)=\sum_{p=0}^{+\infty}a_{2p}t^{2p},$
where\begin{equation}\label{zak7}a_{2p+1}
=-\frac{2(2p-1)\mu}{p\bigl[(2p+1)(2p-1)-(n-1)(n-3)\bigr]}a_{2p-1},\quad\forall\,\,p\geq
1,\end{equation}
\begin{equation}\label{zak8}
\left\{
\begin{array}{ll}
a_{2p}=0, & \qquad\forall\,\,0\leq p< \frac{n-1}{2}; \\
a_{2p}=-\frac{8(p-1)\mu}{(2p-1)\bigl[4p(p-1)-(n-1)(n-3)\bigr]}a_{2p-2},
&\qquad\forall\,\,p\geq \frac{n+1}{2}.
\end{array}
\right.\end{equation} Let $a_1> 0$ and $a_{n-1}< 0$ be fixed. Note
that $\psi_1$ and $\psi_2$ are smooth on $[0,\,1].$ On the other
hand, we claim that, for $\mu$ small enough, we have
\begin{equation}\label{positivity} \psi_1(t)> 0\quad\mathrm{
and }\quad \psi_2(t)< 0,\quad \forall\,\,t\in (0,\,1].\end{equation}
Indeed, it is sufficient to remark that $\psi_1$ and $\psi_2$
satisfy the hypotheses of the alternating series theorem for $\mu$
small enough. Thus there exists a constant $\mu(n)>0$ depending only
on $n,$ such that the claim \eqref{positivity} is valid for every
$\mu\leq \mu(n).$ Denoting
$\bar{\eta}_3:=-2n\bar{c}_3K(y_0)\mu(n)/(n-2)^2\bar{c}_2\Delta
K(y_0).$ \eqref{positivity} enables us to choose $a_1$ and $a_{n-1}$
such that
\begin{equation}\label{zak10}\bar{\psi}(t):=\psi_1(t)+\psi_2(t)\geq 0,\quad \forall \,\,t\in
[0,\,1],\,\,\forall\,\,\mu\leq \mu(n).
\end{equation}
Regarding the identities \eqref{zak5} and \eqref{zak10} and in order
to get the desired contradiction, it is sufficient to investigate a
constant $\bar{\eta}> 0$ such that if $\eta \leq \bar{\eta},$ then,
for any $t\in (0,\,1],$
\begin{equation}\label{zak15}-\frac{n-2}{2n}t\eta f_1'(t)\bar{\psi}(t)+\frac{n-1}{n}
\bigl(f(0)+\eta
f_1(t)\bigr)\Bigl(\bar{\psi}(t)-t\bar{\psi}'(t)\Bigr)>
0.\end{equation} Let $0< \delta\leq 1$ be a fixed constant and
$\delta\leq t\leq 1.$ Combining \eqref{condition}, \eqref{zak7} and
\eqref{zak8} and using the fact that $\mu\leq \mu(n)$ we obtain
\begin{equation}\label{zak19}\begin{split}&\quad-\frac{n-2}{2n}\eta t f_1'(t)\bar{\psi}(t)+\frac{n-1}{n}
f(t)\Bigl(\bar{\psi}(t)-t\bar{\psi}'(t)\Bigr)\\&= \frac{n-1}{n}a_1
f(t)\Bigl(\frac{\eta}{f(0)}O_n(1)-
\frac{(n-2)a_{n-1}}{a_1}t^{n-4}\bigl(1+\frac{\eta}{f(0)}
O_n(1)\bigr)\Bigr)t^3\\&\quad-\frac{n-2}{2n}t\eta
f_1'(t)\bar{\psi}(t),
\end{split}\end{equation}where $|O_n(1)|$ is upper-bounded by a fixed constant $M$ depending only on $n.$ Let $\bar{\eta}_2> 0$ be a constant such that, for any
$0< \eta\leq \bar{\eta}_2,$ \eqref{zak10} is satisfied and
\begin{equation}\label{zak20}-\frac{\eta}{f(0)}\bigl|O_n(1)\bigr|-
\frac{(n-2)a_{n-1}}{a_1}\delta^{n-4}\bigl(1-\frac{\eta}{f(0)}
|O_n(1)|\bigr)> 0.
\end{equation}Combining \eqref{zak10}, \eqref{zak19}, $\mathbf{(K_\eta)},$ and
\eqref{zak20} we obtain \eqref{zak15} for any $\delta\leq t\leq 1$
and any $0< \eta\leq \mathrm{min}(\bar{\eta}_2,\,\bar{\eta}_3).$
Observe that if we let $\delta$ tend to $0,$ then to regain
\eqref{zak20} for $\eta\leq \bar{\eta}_2,$  $\bar{\eta}_2$ must go
to $0:$ this fact leads to the loss of \eqref{zak15}. Thus we have
to fix the constant $\delta$ and we need another argument for the
case $0< t\leq \delta.$ To this end, we will take care of the local
information about the function $f_1(t)$ near its critical point $0.$
First, let us observe that the condition $\mathbf{(K_3)}$ implies
the existence of two constants $\delta,\,M_0> 0$ such that, for any
$0<t\leq \delta,$\begin{equation}\label{M0}0\leq
f_1'(t)-tf_1''(0)\leq M_0t^{n-3}\quad \mathrm{or}\quad
f_1'(t)-tf_1''(0)\leq 0.
\end{equation}In particular, we deduce from \eqref{zak10} and \eqref{M0} that, for any $0<t\leq
\delta,$
\begin{equation}\label{M00}-\bigl(f_1'(t)-tf_1''(0)\bigr)\bar{\psi}(t)\geq-M_0t^{n-3}\bar{\psi}(t)\geq -M_0t^{n-2}(a_1-a_{n-1})|O_n(1)|,
\end{equation}where $|O_n(1)|$ is upper-bounded by a constant
$M_{n}$ depending only on $n.$ Now, by combining \eqref{condition},
\eqref{zak7}, \eqref{zak8} and \eqref{zak10} and using the fact that
$\mu\leq \mu(n)$ we
obtain\begin{equation}\label{zak16}\begin{split}&-\frac{n-2}{2n}\eta
tf_1'(t)\bar{\psi}(t)+\frac{n-1}{n} \bigl(f(0)+\eta
f_1(t)\bigr)\Bigl(\bar{\psi}(t)-t\bar{\psi}'(t)\Bigr)\\&=\Bigl(-\frac{n-2}{2n}\eta
f_1''(0)a_1-2f(0)\frac{(n-1)}{n}a_3\Bigr)t^3-\frac{n-2}{2n}\eta
t\bigl(f_1'(t)-tf_1''(0)\bigr)\bar{\psi}(t)\\&+\Bigl[-\bigl(\frac{n-2}{2n}+\frac{n-1}{n}\bigr)\eta
f_1''(0) a_3-4f(0)\frac{n-1}{n}a_5+\frac{\eta^2\mu}{f(0)}a_1
O_{n,f_1}(1)\Bigr]t^5\\&\,\,+f(0)\Bigl(-\frac{(n-2)(n-1)}{n}
a_{n-1}+\frac{\eta}{f(0)} O_{n,f_1}(1)(-a_{n-1}+a_1)\Bigr)t^{n-1},
\end{split}\end{equation}
where $|O_{n,f_1}(1)|$ is upper-bounded by a fixed constant
$M_{n,f_1}$ depending only on $n$ and the function $f_1(x).$
Finally, by using \eqref{zak7}, \eqref{condition} and the fact that
$n\neq 7\text{--}19$ and that
$\bar{c}_2/\bar{c}_3=n(n-4)/4(n-1)(n-2)$ we get
\begin{equation}\label{zak17}
-\frac{n-2}{2n}\eta f_1''(0) a_1-2f(0)\frac{n-1}{n}a_3\geq 0,
\end{equation}
\begin{equation}\label{zak18}\frac{1}{\eta\mu a_1}\Bigl(-\bigl(\frac{n-2}{2n}+\frac{n-1}{n}\bigr)\eta f_1''(0)
a_3-4f(0)\frac{n-1}{n}a_5\Bigr)\geq M> 0,
\end{equation}where $M$ is a constant depending only on
$n.$ Combining \eqref{M00}--\eqref{zak18} and taking $\bar{\eta}_1>
0$ small enough such that, for any $0< \eta\leq
\bar{\eta}_1,$\begin{eqnarray}
-\frac{(n-2)(n-1)}{n} a_{n-1}-\frac{\eta}{f(0)}
(-a_{n-1}+a_1)\bigl(|O_{n,f_1}(1)|+\frac{n-2}{2n}M_0|O_n(1)|\bigr)>
0,\nonumber\\ M-\frac{\eta}{f(0)}|O_{n,f_1}(1)|>
0,\qquad\qquad\qquad\qquad\qquad\qquad\qquad\qquad\qquad\qquad\qquad\quad\,\,\,\,\nonumber
\end{eqnarray}we get \eqref{zak15} for any $0< t\leq \delta$ and
any $0< \eta\leq \mathrm{min}(\bar{\eta}_1,\,\bar{\eta}_3).$ The
proof of Theorem \ref{nonexistence} follows by choosing
$\bar{\eta}=\mathrm{min}(\bar{\eta}_1,\,\bar{\eta}_2,\,\bar{\eta}_3).$\\\\
{\bf Acknowledgement.} The author is greatly indebted to Professor
A. Bahri, Professor H. Brezis, Professor P. L. Lions, and Professor
L. Nirenberg for their works that gave him the inspiration to
prepare this work.


\begin{thebibliography}{99}
\bibitem{A}
T. Aubin, {\it Equations diff\'erentielles non lineaires et
probl\`eme de Yamabe concernant la courbure scalaire}, J. Math.
Pures Appl. {\bf 55} (1976), no. 3, 269--293.

\bibitem{B}
A. Bahri, {\it Critical points at infinity in some variational
problems}, Pitman Research Notes in Mathematics Series, Vol. 182,
Longman Scientific and Technical, Harlow, 1989.

\bibitem{Bo}
Z. Boucheche, {\it Existence result for an elliptic equation
involving critical exponent in three-dimensional domains}, Complex
Var. Elliptic Equ. {\bf 64} (2019), no. 4, 649--675.

\bibitem{Bo1}
Z. Boucheche, {\it .........}, Submitted.

\bibitem{BrNir}
H. Brezis and L. Nirenberg, {\it Positive solutions of nonlinear
elliptic equations involving critical Sobolev exponents}, Comm. Pure
Appl. Math. {\bf 36} (1983), no. 4, 437--477.

\bibitem{GNNir}
B. Gidas, W. M. Ni and L. Nirenberg, {\it Symmetry and related
properties via the maximum principle}, Comm. Math. Phys. {\bf 68}
(1979), no. 3, 209--243.


\bibitem{Li}
P. Lions, {\it The concentration compactness principle in the
calculus of variations. The limit case, part 2}, Rev. Mat. Iberoam.
{\bf 1} (1985), no. 2, 45--121.

\bibitem{P}
S. Pohozaev, {\it Eigenfunctions of the equation $\Delta u+\lambda
f(u)=0$}, Soviet Math. Doklady, {\bf 6} (1965), 1408--1411.
\end{thebibliography}
\end{document}